# Generation of acyclic biological diagrams

**A. Panayotopoulos**

**Abstract** For the generation of acyclic biological diagrams, from a graph-theoretical perspective, we introduce the relative diagrams of cyclic permutations with ramphoid and keratoid vertices of degree 2, which correspond to Motzkin and Dyck words/paths. The relation between these two types of diagrams, defines the generation of the first via the permutations of the second, which may be of assistance in the study and treatment of several biological problems.

## 1 Introduction

One of the most important research tools in Biology nowadays are the diagrams, which may take multiple forms, according to their biological context [5] [6] [7]. This paper examines the acyclic biological diagrams ($b$-diagrams for short - see Fig. 8) which are graphs over the vertex set $[n]=\{1, 2, \ldots, n\}$ with a vertex degree at most 2 and are drawn as points on a horizontal line, while their edges are drawn as arcs on the upper half plane [6].

In Section 2, the diagrams of the cyclic permutations $\sigma$ ($\sigma$-diagrams for short - see Fig. 1) are graphically represented by Motzkin and Dyck paths. The solution to the problem of finding the subset of cyclic permutations which correspond to a related Motzkin or Dyck word, is also examined.

Section 3 examines the $b$-diagrams which have common axis and the same number of vertices as the $\sigma$-diagrams, but a vertex degree of at most 2 and fewer arcs. These $b$-diagrams are graphically represented with paths that look like Motzkin and Dyck paths, without being such, since their alphabet uses more letters. The solution to the problem of finding the $b$-diagrams that correspond to related words is provided. The algebra of $b$-diagrams is presented and the notions of kernel, inflation and the $k$-noncrossing structures are also introduced.

Section 4 deals with the problem of finding the generators of a $b$-diagram, i.e., of the cyclic permutations for which the arc set of their diagrams is a superset of the arc set of the respective $b$-diagrams. The number of these generators is calculated and the notion of the complement of every $b$-diagram is defined. A reference is made to Plato's transformations and finally, in the conclusions, the possible applications in biological cases, are suggested.

## 2 Diagrams of cyclic permutations

We denote by $\sigma = \sigma_1\sigma_2 \ldots \sigma_n$, where $\sigma_1=1$, every cyclic permutation of the set $[n] = \{1, 2, \ldots, n\}$. The set of all cyclic permutations of $[n]$ is denoted by $\Sigma_n$, where $|\Sigma_n|=(n-1)!$ . For convenience we extend the set of indices from $[n]$ to the set of all integers, by defining $\sigma_i = \sigma_j$ whenever $i \equiv j \pmod n$, e.g., $\sigma_n = \sigma_0$ and $\sigma_{n+1} = \sigma_1$.

Every $\sigma \in \Sigma_n$ corresponds to a cyclic graph with vertex set $[n]$ and edge set $\{\{\sigma_1, \sigma_2\}, \{\sigma_2, \sigma_3\},\ldots,\{\sigma_{n-1}, \sigma_n\}, \{\sigma_n, \sigma_1\}\}$

A. Panayotopoulos ( 🖂 )
University of Piraeus, Karaoli & Dimitriou 80, 18534, Piraeus, Greece
e-mail: antonios@unipi.gr

drawn as in Fig. 1(a). Obviously, these edges define a Hamiltonian circuit. If we illustrate the previous graph with its vertices drawn as points on a horizontal line and its edges drawn as arcs above this line, the resulting representation is called the $\sigma$-diagram of $\sigma$ (Fig. 1(b)).
The reverse of $\sigma$ is the cyclic permutation $\bar{\sigma} = \sigma_1\sigma_n\sigma_{n-1}\sigma_2$. It can be obtained by setting $\bar{\sigma}_i = \sigma_{n+2-i}$, $i \in [n]$ and it has the same $\sigma$-diagram with $\sigma$. For example, the reverse of $\sigma = 1\,3\,2\,7\,8\,4\,5\,6$ is $\bar{\sigma} = 1\,6\,5\,4\,8\,7\,2\,3$.

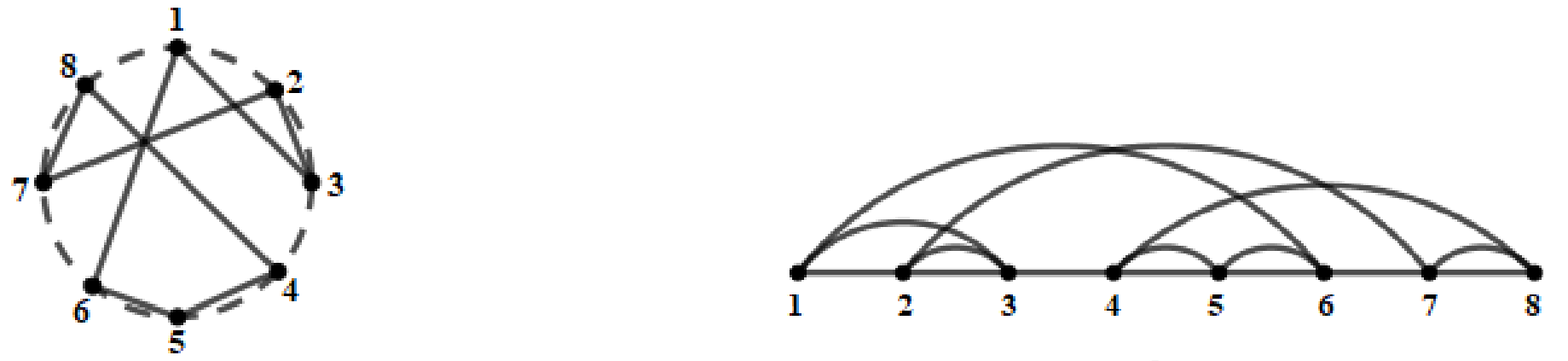

(a) (b)

Figure 1: The cyclic graph (a) and the $\sigma$-diagram of $\sigma = 1\,3\,2\,7\,8\,4\,5\,6$ (b)

Clearly, the arcs of a $\sigma$-diagram form a cycle and therefore each vertex has degree 2. Note, that each $\sigma$-diagram corresponds to exactly two cyclic permutations, each one being the reverse of the other. Moreover, a $\sigma$-diagram is equivalently defined by the set $U_\sigma$ of its ordered pairs of arcs (i.e. the arc corresponding to one edge $\{i, j\}$ with $i < j$, is denoted by $(i, j)$ or $ij$ for simplicity). In the $\sigma$-diagram of Fig. 1(b), we have $U_\sigma = \{13, 16, 23, 27, 45, 48, 56, 78\}$. Obviously, $U_{\bar{\sigma}} = U_\sigma$ .
It is easy to see that each vertex of an $\sigma$-diagram belongs to exactly one of the following sets of classes [3]:

- The set of left ramphoids $R = \{\sigma_i : \sigma_i < \min\{\sigma_{i-1}, \sigma_{i+1}\},\ i \in [n]\}$.
- The set of right ramphoids $\bar{R} = \{\sigma_i : \sigma_i < \max\{\sigma_{i-1}, \sigma_{i+1}\},\ i \in [n]\}$.
- The set of keratoids $K = \{\sigma_i : \sigma_{i-1} < \sigma_i < \sigma_{i+1} \text{ or } \sigma_{i-1} > \sigma_i > \sigma_{i+1},\ i \in [n]\}$.

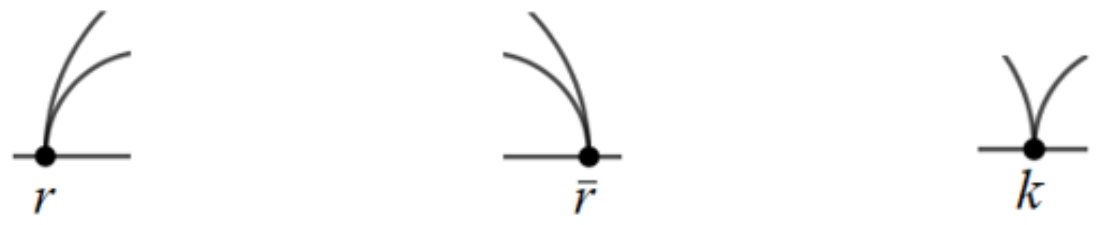

Figure 2: The 3 types of vertices: left ramphoid ($r$), right ramphoid ($\bar{r}$) and keratoid ($k$).

For example, in the $\sigma$-diagram of Fig. 1(b), we have $R = \{1, 2, 4\}$, $\bar{R} = \{3, 6, 8\}$ and $K = \{5, 7\}$, according to the inequalities $(6 >)\ 1 < 3 > 2 < 7 < 8 > 4 < 5 < 6\ (> 1)$. Similarly, for the $\sigma$-diagram of Fig. 3, we have $R = \{1,2,4,5\}$, $\bar{R} = \{3,6,7,8\}$ and $K = \emptyset$.
It must be emphasized that for the purposes of this paper, the arcs of the ramphoid vertices will always form nestings and never crossings.

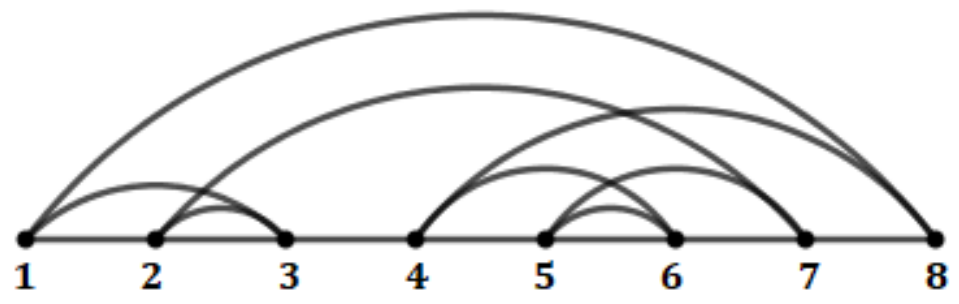

Figure 3: The $\sigma$-diagram of $\sigma = 1\,3\,2\,7\,5\,6\,4\,8$

The sets $R$, $\bar{R}$ and $K$ can uniquely be determined by using only the set of arcs $U_\sigma$ of the $\sigma$-diagram. The elements of $R$ (resp. $\bar{R}$) are the vertices that appear twice in $U_\sigma$ as the starting vertex $i$ (resp. as the ending vertex $j$) of an arc $(i,j)$. For example, in the $\sigma$-diagram of Fig. 1(b), we have:

| | | | | | | | | |
|---|---|---|---|---|---|---|---|---|
| $i$ | 1 | 1 | 2 | 2 | 4 | 4 | 5 | 7 |
| $j$ | 3 | 6 | 3 | 7 | 5 | 8 | 6 | 8 |

where $R$ (resp. $\bar{R}$) consists of the elements appearing twice in the first (resp. second) row, i.e., $R = \{1,2,4\}$ and $\bar{R} = \{3,6,8\}$, whereas $K$ consists of the remaining elements, i.e. $K = \{5,7\}$. It is easy to prove that for every $\sigma$-diagram $|R| = |\bar{R}|$. Since the elements of $R$ (resp. $\bar{R}$) are common starting points (resp. ending points) of arc pairs and the elements of $K$ are starting and ending points of consecutive arcs, it is: $|R| + \left|\frac{K}{2}\right| = |\bar{R}| + \left|\frac{K}{2}\right|$. So, since $R \cup \bar{R} \cup K = [n]$, we have $2|R| + K = n$ which means that $n$ and $|K|$ have the same parity.

**Motzkin $\sigma$-diagrams**

Given an $\sigma$-diagram with $K \neq \emptyset$, if we assign the letters $r$, $\bar{r}$ and $k$ to the elements of $R$, $\bar{R}$ and $K$ respectively, then we obtain a word $w = w_1 w_2 \ldots w_n \in \{r, \bar{r}, k\}^*$, such that for every $i \in [n]$:

$w_i = \begin{cases} r, & \text{if } i \in R \\ \bar{r}, & \text{if } i \in \bar{R} \\ k, & \text{if } i \in K \end{cases}$. For example, the $\sigma$-diagram of Fig. 1(b) corresponds to the word $w = r r \bar{r} r k \bar{r} k \bar{r}$.

The set of all those words $w \in \{r, \bar{r}, k\}^*$ is denoted by $W_n$, to which the following rules apply:

i) For every $w \in W_n$, we have $w_1 = r$ and $w_2 \neq \bar{r}$ (resp. $w_{n-1} \neq r$ and $w_n = \bar{r}$).

ii) Every $w \in W_n$, has two corresponding sets of permutations $\sigma \in \Sigma_n$. The permutations that have $\sigma_2 = \min \bar{R} \cup K$ belong to the first set, whereas their reverse ones belong to the second set.

Recall that a word $w \in \{r, \bar{r}, k\}^*$ is a Motzkin word iff $|w|_r = |w|_{\bar{r}}$ and $|p|_r \geq |p|_{\bar{r}}$ for every factorization $w = pq$, where $|w|_x$ denotes the number of occurrences of the letter $x$ in the word $w$ and $|w|$ denotes the length of $w$. By assigning the letter $r$ to an up-step (1,1), $\bar{r}$ to a down-step (1,-1) and $k$ to a flat step (1,0), a Motzkin word corresponds to a unique Motzkin path, starting at the origin (0,0) and ending at $(n, 0)$ on the $x$-axis, (see Fig. 4). It is called Motzkin $\sigma$-diagram and it's path touches the $x$-axis only at its two endpoints, so it is an elevated Motzkin path, hence $|W_n| = M_{n-2}$, where $M_n$ is the **n**-th motzkin number (sequence A001006 in [8]).

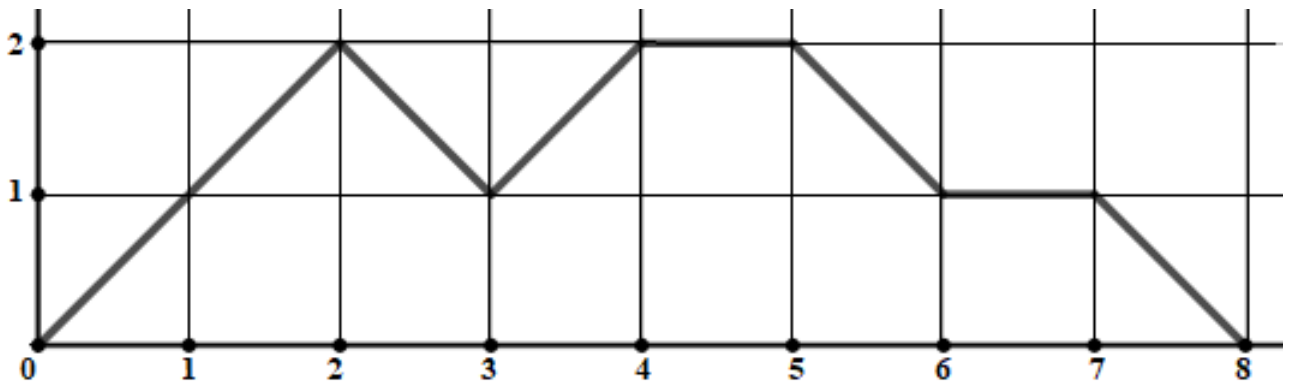

Figure 4: The Motzkin path corresponding to the word $w = r\, r\, \bar{r}\, r\, k\, \bar{r}\, k\, \bar{r}$.

iii) If $w' = w'_1 w'_2 \ldots w'_n$ where $w'_i = w_{\sigma_i}$, are determined not by the natural order, but according to the order indicated by $\sigma$, then $w_i = w'_{\sigma_i^{-1}}$, where $\sigma^{-1}$ denotes the inverse of $\sigma$, i.e., $\sigma_j^{-1} = i \Leftrightarrow \sigma_i = j$. Indeed, by using the identity $\sigma(\sigma_i^{-1}) = i$, we have $w'_{\sigma_i^{-1}} = w_{\sigma(\sigma_i^{-1})} = w_i$.

| $i$ | 1 | 2 | 3 | 4 | 5 | 6 | 7 | 8 |
|---|---|---|---|---|---|---|---|---|
| $\sigma_i$ | 1 | 3 | 2 | 7 | 8 | 4 | 5 | 6 |
| $\sigma_i^{-1}$ | 1 | 3 | 2 | 6 | 7 | 8 | 4 | 5 |
| $w'_i$ | $r$ | $\bar{r}$ | $r$ | $k$ | $\bar{r}$ | $r$ | $k$ | $\bar{r}$ |
| $w_i$ | $r$ | $r$ | $\bar{r}$ | $r$ | $k$ | $\bar{r}$ | $k$ | $\bar{r}$ |

**Dyck $\sigma$-diagrams**

Given an $\sigma$-diagram with $K = \emptyset$ (see Fig. 3), then we obtain a word $u = u_1 u_2 \ldots u_n$ of length $n = 2|R|$ such that for every $i \in [n]$ we have $u_i = \begin{cases} r, & i \in R \\ \bar{r}, & i \in \bar{R} \end{cases}$.

For example, the $\sigma$-diagram of Fig. 3 corresponds to the word $u = r r \bar{r} r r \bar{r} \bar{r} \bar{r}$.

The set of all the above words $u \in \{r, \bar{r}\}^*$ is denoted by $U_n$, for which we have:

i) For every $u \in U_n$ we have $u_1 = r$ and $u_n = \bar{r}$.

ii) Each $u \in U_n$ is a Dyck word of length $n$ (see Fig. 5).

It is obvious that the number $|U_n|$ coinsides with the number of all the Dyck paths (sequence A057163 in [8]).

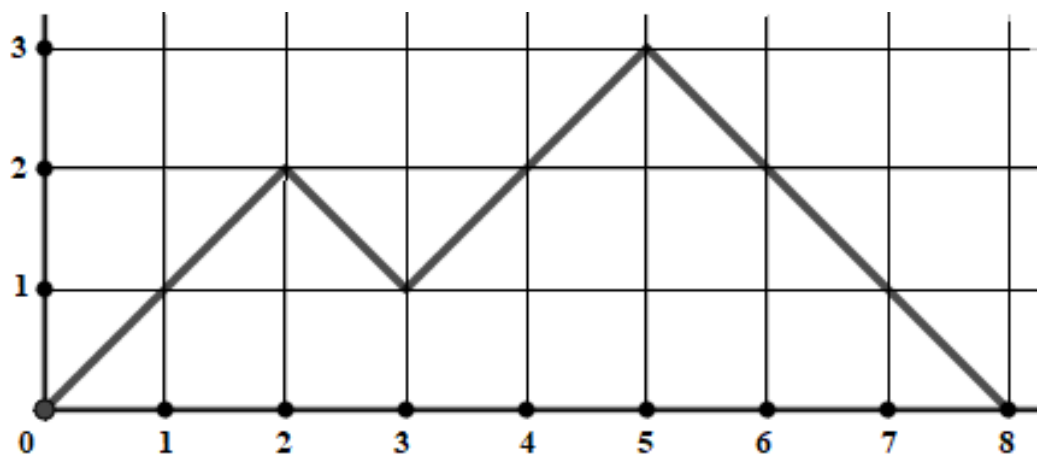

Figure 5: The Dyck path corresponding to the word $u = r\,r\,\bar{r}\,r\,r\,\bar{r}\,\bar{r}\,\bar{r}$.

iii) For every $u \in U_n$ it is: $u_i = \begin{cases} r, & \text{if } \sigma_i^{-1} \text{ is odd} \\ \bar{r}, & \text{if } \sigma_i^{-1} \text{ is even} \end{cases}$ because, for the corresponding word $u' = u'_1 u'_2 \ldots u'_n$ with $u'_i = u_{\sigma_i}$. Thus, the following table derives for the example of Fig. 5:

| $i$ | 1 | 2 | 3 | 4 | 5 | 6 | 7 | 8 |
|---|---|---|---|---|---|---|---|---|
| $\sigma_i$ | 1 | 3 | 2 | 8 | 4 | 7 | 5 | 6 |
| $\sigma_i^{-1}$ | 1 | 3 | 2 | 5 | 7 | 8 | 6 | 4 |
| $u'_i$ | $r$ | $\bar{r}$ | $r$ | $\bar{r}$ | $r$ | $\bar{r}$ | $r$ | $\bar{r}$ |
| $u_i$ | $r$ | $r$ | $\bar{r}$ | $r$ | $r$ | $\bar{r}$ | $\bar{r}$ | $\bar{r}$ |

**Problem of words**

Every permutation $\sigma \in \Sigma_n$ corresponds to a unique word $w \in W_n$ (resp. $u \in U_n$), but the opposite is not valid, i.e. the mapping is not injective. Thus we have the problem of finding the subset of $\Sigma_n$ which corresponds to a word $w \in W_n$ (resp. $u \in U_n$).

i) For every vertex $i \in [n]$ of a Motzkin $\sigma$-diagram we have the following correspondences:

$$\Gamma(i) = \begin{cases} \{ j \in \bar{R} \cup K \text{ and } j > i \} & \text{if } i \in R \\ \{ j \in R \cup K \text{ and } j < i \} & \text{if } i \in \bar{R} \\ \{ j \in R \cup K \text{ and } j < i \} \cup \{ j \in \bar{R} \cup K \text{ and } j > i \} & \text{if } i \in K \end{cases}$$

So, for the word $w = r\,k\,r\,\bar{r}\,k\,\bar{r}$ of $W_6$ it is $R = \{1,3\}$, $\bar{R} = \{4,6\}$, $K = \{2,5\}$, $R \cup K = \{1,2,3,5\}$, $\bar{R} \cup K = \{2,4,5,6\}$ and $\Gamma(1)=\{2,4,5,6\}$, $\Gamma(2)=\{1,4,5,6\}$, $\Gamma(3)=\{4,5,6\}$, $\Gamma(4)=\{1,2,3\}$, $\Gamma(5)=\{1,2,3,6\}$, $\Gamma(6)=\{1,2,3,5\}$. From the tree of Fig. 6 we get the permutations 1 2 4 3 5 6, 1 2 4 3 6 5, 1 2 5 6 3 4 and 1 2 6 5 3 4 which make up half of the total number of permutations, whereas the other half are their reverse ones: 1 6 5 3 4 2, 1 5 6 3 4 2, 1 4 3 6 2 5 and 1 4 3 5 6 2.

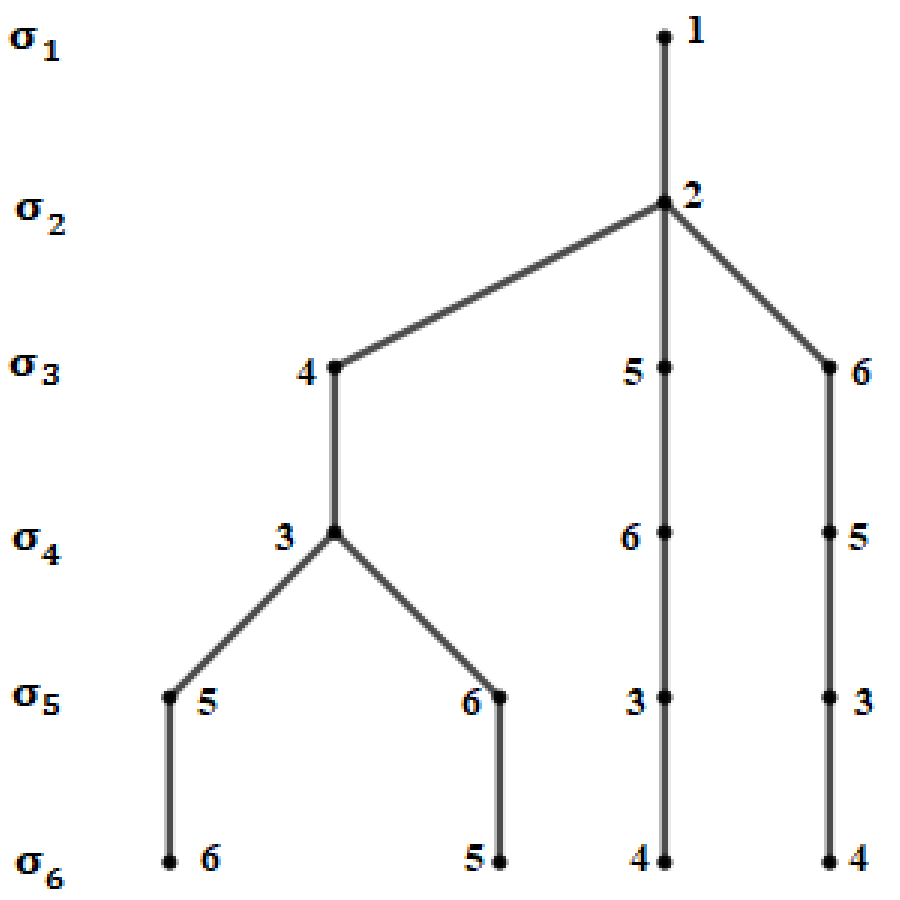

Figure 6: The tree of permutations 1 2 4 3 5 6, 1 2 4 3 6 5, 1 2 5 6 3 4 and 1 2 6 5 3 4.

ii) Similarly, for every vertex $i \in [n]$ of a Dyck $\sigma$-diagram we have the following correspondences:

$$\Gamma(i) = \begin{cases} \{ j \in \bar{R} \text{ and } j > i \} & \text{if } i \in R \\ \{ j \in R \text{ and } j < i \} & \text{if } i \in \bar{R} \end{cases}$$

So, for the word $u = r\,r\,\bar{r}\,r\,\bar{r}\,\bar{r}$ of $U_6$ it is $R = \{1,2,4\}$, $\bar{R} = \{3,5,6\}$ and $\Gamma(1)=\{3,5,6\}$, $\Gamma(2)=\{3,5,6\}$, $\Gamma(3)=\{1,2\}$, $\Gamma(4)=\{5,6\}$, $\Gamma(5)=\{1,2,4\}$, $\Gamma(6)=\{1,2,4\}$. From the tree of Fig. 7 we get the permutations 1 3 2 5 4 6 and 1 3 2 6 4 5, which also make up half of the total number of permutations, whereas the other half are their reverse ones: 1 6 4 5 2 3 and 1 5 4 6 2 3.

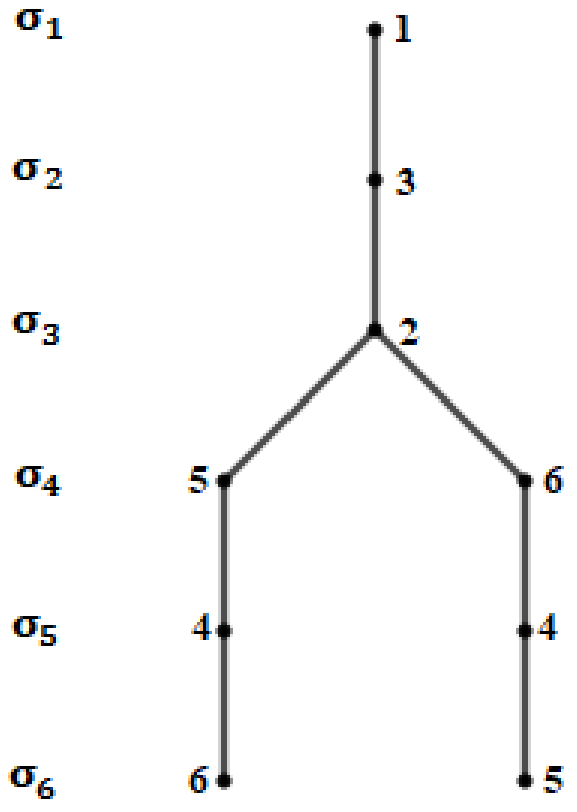

Figure 7: The tree of permutations 1 3 2 5 4 6 and 1 3 2 6 4 5.

## 3 Biological Diagrams

In Biology, several models (secondary structures, braids, tangles, pseudoknots etc.) can be represented and studied via diagrams similar to $\sigma$-diagrams. These diagrams are described by specific definitions about the nature and the role of the horizontal line, the vertices and the arcs, according to their biological context. From a theoretical perspective, they are acyclic subgraphs of $\sigma$-diagrams and they derive from them after the removal of a set of arcs $C$. For example, if from the $\sigma$-diagram of $\sigma = 1\ 3\ 2\ 7\ 8\ 4\ 5\ 6$ (Fig. 1(b)) we remove the set of arcs $C = \{27,48,56\}$ (resp. $C = \{16, 23, 27, 45, 78\}$), the diagrams of Fig. 8(a) $b = 3\ 1\ 6 | 2\ 7\ 8 | 4\ 5$ (resp. of Fig. 8(b) $b = 1\ 3 | 2 | 4\ 8 | 5\ 6 | 7$ ) derive.

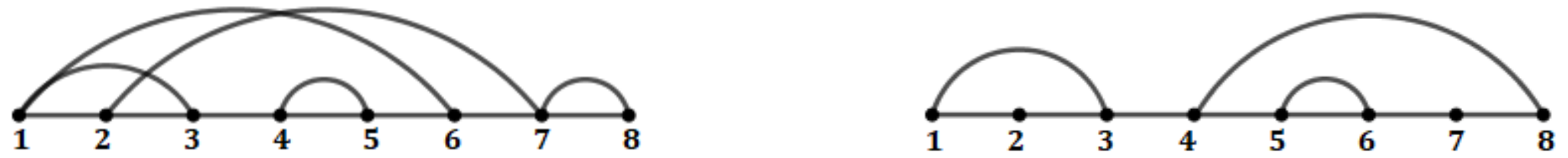

(a) (b)

Figure 8: $b = 3\ 1\ 6 | 2\ 7\ 8 | 4\ 5$ $b = 1\ 3 | 2 | 4\ 8 | 5\ 6 | 7$

In general, a $b$-diagram in $[n]$ is a partition of any permutation of $[n]$ into $m$ subsets $b_1, b_2, \ldots, b_m$, each one of which includes $n_i \in [n\text{ -}1]$, $i \in [m]$ consecutive terms (when $n_i > 1$) or isolated terms (when $n_i = 1$), where $n_1 + n_2 + \ldots + n_m = n$. This diagram will be denoted by $b = b_1 | b_2 | \ldots | b_m$, where $b_i$, $i \in [m]$ are called its blocks and correspond to its arcs (when $n_i$=2), its paths (when $n_i \geq 3$), and its isolated vertices [(1)] (when $n_i$=1), while in every block correspond $n_i - 1$ consecutive arcs. It is easy to derive that each $b$-diagram can be fully defined by the set $U_b$ of its ordered arcs and its isolated vertices. Thus, the set $U_b$={13,16,27,45,78}, fully defines the $b$-diagram of Fig. 8(a) which forms a braid with $C$ ={23,48,56}. Respectively, the set $U_b$={13,2,48,56,7}, fully defines the $b$-diagram of Fig. 8(b) which forms a secondary structure with $C$={16,23,27,45,78}. The set-difference $C = U_\sigma \setminus U_b$ which is called Cut Set, characterizes the relation between the $\sigma$-diagrams and the $b$-diagrams, where $|U_\sigma| = n$ and $|U_b| = n - m$. The set of vertices of a $b$-diagram is partitioned into 6 classes, 3 of which are $R$, $\bar{R}$ and $K$. The remaining 3 classes are the class $L$ of isolated vertices, the class $A$ of vertices which are beginnings of unique arcs and the class $\bar{A}$ of vertices which are ends of unique arcs. By assigning to these vertices the letters $e$, $\alpha$ and $\bar{\alpha}$ respectively, each $b$-diagram corresponds to a word $z = z_1 z_2 \ldots z_n$ with $z_i \in \{e, \alpha, \bar{\alpha}, r, \bar{r}, k\}$.

---

*(1): which derive from the removal of consecutive arcs from the $\sigma$-diagram.*

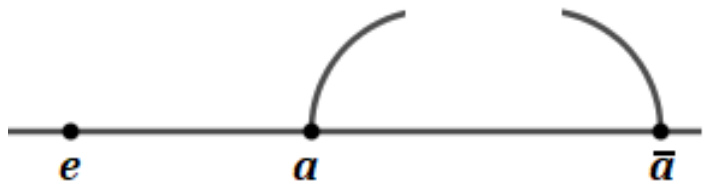

Figure 9: The 3 types of vertices in $L$, $A$ and $\bar{A}$.

For example, the words $z = r\,\alpha\,\bar{\alpha}\,\alpha\,\bar{\alpha}\,\bar{\alpha}\,k\,\bar{\alpha}$ and $z = \alpha\,e\,\bar{\alpha}\,\alpha\,\alpha\,\bar{\alpha}\,e\,\bar{\alpha}$ correspond to the $b$-diagrams of Fig. 8(a) and Fig. 8(b), respectively. The set of the above words $z \in \{e, \alpha, \bar{\alpha}, r, \bar{r}, k\}^*$ which correspond to $b$-diagrams will be symbolized with $Z_n$, therefore in general it will be $z_1 \in \{\alpha, r\}$ and $z_n \in \{\bar{\alpha}, \bar{r}\}$.

Property: For every $b$-diagram it is $|A|+2|R| = |\bar{A}|+2|\bar{R}|$

Property: The length of every word $z \in Z_n$ is $n+|R|+|\bar{R}|+|K|$.

To every vertex $i \in [n]$ of a $b$-diagram with $z = z_1 z_2 \ldots z_n$, corresponds a number $\theta_i \in \{-2, -1, 0, 1, 2\}$ which is called degree, such as:

$$\theta_i = \begin{cases} -2 & \text{if } z_i = \bar{r} \\ -1 & \text{if } z_i = \bar{\alpha} \\ 0 & \text{if } z_i = e \text{ or } k \\ 1 & \text{if } z_i = \alpha \\ 2 & \text{if } z_i = r \end{cases}$$

Thus, we get the vector $\theta = (\theta_1, \theta_2, \ldots, \theta_n)$ of the degrees of the vertices.

Proposition: For every $b$-diagram it is $\sum_{i=1}^{n} \theta_i = 0$.

The above is valid for every $\sigma$-diagram with $\theta_i \in \{-2, 0, 2\}$ because of the equation $|R|=|\bar{R}|$ and it will also be valid for every $b$-diagram which derives from these $\sigma$-diagrams after the removal of some arcs, i.e. pairs of $(+1, -1)$ that produce a zero sum. Subsequently, the following necessary and sufficient condition is also valid:

Proposition: A word $z \in Z_n$ defines a $b$-diagram iff $\sum_{i=1}^{k} \theta_i \geq 0$ for every $k \leq n-1$ and $\sum_{i=1}^{n} \theta_i = 0$.

Thus, a $b$-diagram exists for the word $z = r\,a\,r\,\bar{a}\,\bar{r}\,\bar{a}\,\bar{a}$, but not for the word $\zeta = \bar{r}\,\bar{a}\,k\,e\,\alpha\,r$, although $\sum_{i=1}^{6} \theta_i = 0$.

Every $b$-diagram has a corresponding path which is not a Motzkin nor a Dyck path, since the alphabet has 6 letters instead of 3. The letters of every word $z = z_1 z_2 \ldots z_n$ in $\{e, \alpha, \bar{\alpha}, r, \bar{r}, k\}^*$ are represented as in Fig. 10.

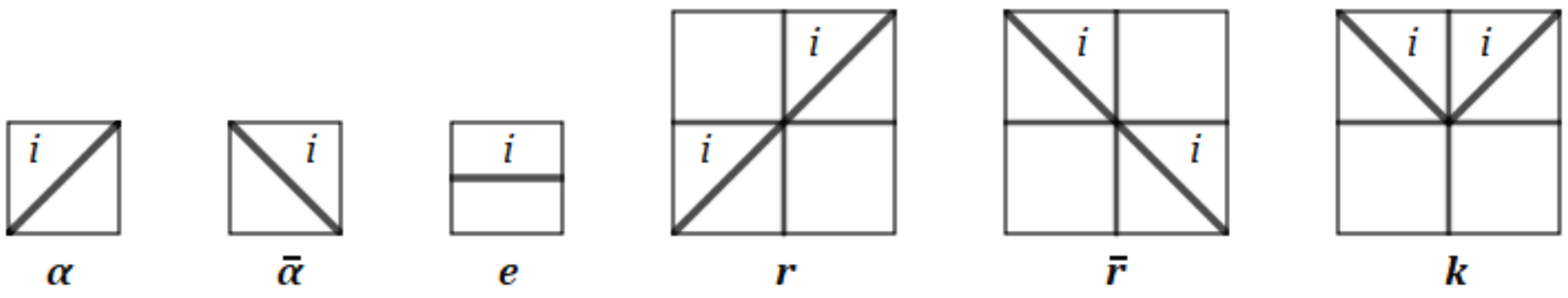

Figure 10: The representations of $\alpha, \bar{\alpha}, e, r, \bar{r}, k$:

$\alpha$: single ascending diagonal line,
$\bar{\alpha}$: single descending diagonal line,
$e$: single horizontal line,
$r$: double ascending diagonal line,
$\bar{r}$: double descending diagonal line,
$k$: single descending diagonal line followed by a single ascending diagonal line.

In all the above representations, the index $i$ of the respective vertex, is also noted.

(a) 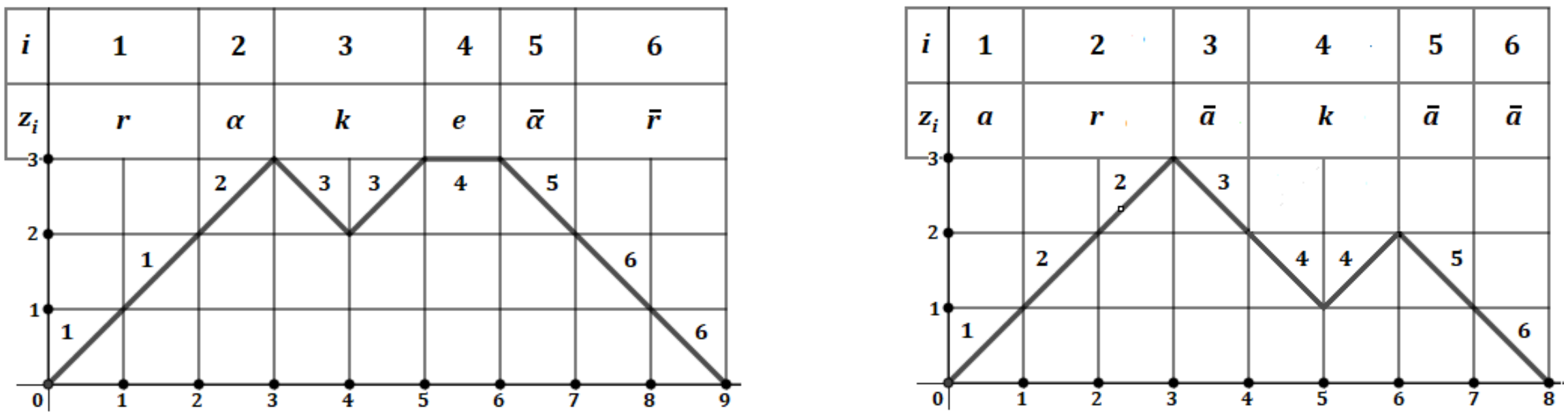 (b)

Figure 11: Paths of $z = r\,\alpha\,k\,e\,\bar{\alpha}\,\bar{r}$ (a) and $z = \alpha\,r\,\bar{\alpha}\,k\,\bar{\alpha}\,\bar{\alpha}$ (b)

The $b$-diagrams associated with a word $z \in \{\alpha, \bar{\alpha}, r, \bar{r}, k\}^*$ can be determined by the following correspondences:

$$\Gamma(i) = \begin{cases} \{ j \in \bar{A} \cup \bar{R} \cup K \text{ and } j > i \} \text{ if } i \in A \cup R \\ \{ j \in A \cup R \cup K \text{ and } j < i \} \text{ if } i \in \bar{A} \cup \bar{R} \\ \{ j \in A \cup R \cup K \text{ and } j < i \} \cup \{ j \in \bar{A} \cup \bar{R} \cup K \text{ and } j > i \} \text{ if } i \in K \end{cases}$$

Thus, for $z = \alpha\, r\, \bar{\alpha}\, k\, \bar{\alpha}\, \bar{\alpha}$ it is $A = \{1\}$, $\bar{A} = \{3,5,6\}$, $R = \{2\}$, $\bar{R} = \emptyset$, $K = \{4\}$, $A \cup R \cup K = \{1,2,4\}$, $\bar{A} \cup \bar{R} \cup K = \{3,5,6,4\}$ and $\Gamma(1)=\{3,4,5,6\}$, $\Gamma(2)=\{3,4,5,6\}$, $\Gamma(3)=\{1,2\}$, $\Gamma(4)=\{1,2,5,6\}$, $\Gamma(5)=\{1,2,4\}$, $\Gamma(6)=\{1,2,4\}$. The paths of the tree of Fig. 12 define the respective $b$-diagrams, whereas the dotted line separates their blocks. The paths that do not produce a $b$-diagram have been omitted (for example the paths that have 4 as the start and the end of a block).

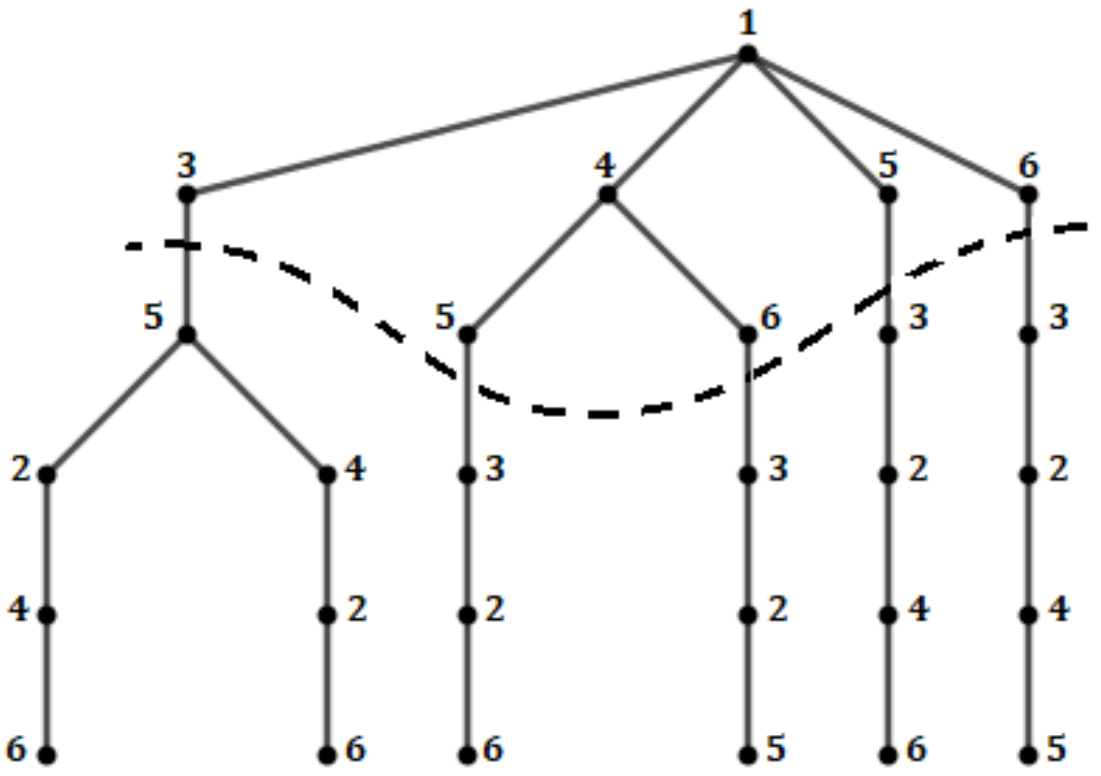

Figure 12: Tree of the $b$-diagrams 13|5246 , 13|5426 , 145|326 , 146|325, 15|3246, 16|3245

The algebra of $b$-diagrams is defined by their arcs set $U_b$ where $U_b \subset U_\sigma$ . If $b$ and $b'$ are two $b$-diagrams in $[n]$ and $[n']$ respectively, it is:

i) Union: $b \cup b'$ we call the $b$-diagram with $U_{b \cup b'} = U_b \cup U_{b'}$ in [max{n,$n'$}]

ii) Intersection: $b \cap b'$ we call the $b$-diagram with $U_{b \cap b'} = U_b \cap U_{b'}$ in [min{n,$n'$}]

iii) When $b \cap b' = \emptyset$ then their $b$-diagrams are disjoint.

iv) Concatenation: $b\,b'$ we call the $b$-diagram with $U_{bb'} = U_b \cup U_{b^*}$ where $b^*$ is the $b'$ in $[n+1, n+2, \ldots, n+n']$.

**Kernel**

Kernel [1] of a $b$-diagram in $[n]$ is called every set $\pi \subset [n]$ which is stable (i.e. its vertices are not connected with arcs) and absorbent (i.e. its vertices are connected with arcs with the vertices of $[n] \backslash \pi$).
Thus, for the $b$-diagrams $b$=1 4 | 2 5 6 | 3 8 7 (Fig 13(a)) and $b$=1 6 3 4 5 | 2 7 10 8 | 9 11 12 (Fig 13(b)), their respective Kernels are $\pi = \{1,5,8\}$ and $\pi = \{4,6,7,8,11\}$.

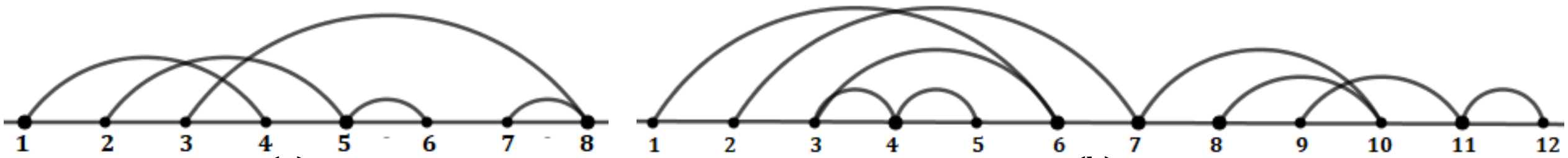

(a) (b)

Figure 13: $b$= **1** 4 | 2 **5** 6 | 3 **8** 7   $b$= 1 **6** 3 **4** 5 | 2 **7** 10 **8** | 9 **11** 12

It must be noted that a kernel is not necessarily unique for a $b$-diagram. Thus, for the $b$-diagram of Fig 13(a), the set $\pi = \{1,2,6,8\}$ is also a kernel. The kernels with $\min|\pi|$ are also of interest and if there exist any isolated vertices, they will belong to a kernel. Finally, the following are valid: $|\pi| \geq n\text{-}m$ and $\pi \subset A$.

**Inflation**

Given a certain $b$-diagram of a word $z \in Z_n$ without isolated vertices, then the $b$-diagram which derives from it by applying the transformations $r \to \alpha\,\alpha$, $\bar{r} \to \bar{\alpha}\,\bar{\alpha}$ and $k \to \bar{\alpha}\,\alpha$, is called inflation (see [6] p. 68). It is obvious that to every inflation will correspond a Dyck word $a \in \{\alpha, \bar{\alpha}\}^*$ with length $|a| = n + |R| + |\bar{R}| + |K|$. Thus, to the $b$-diagram of $z = \alpha\, r\, \bar{\alpha}\, k\, \bar{\alpha}\, \bar{\alpha}$ (Fig. 14(a)), corresponds as its inflation the $b$-diagram of $a = \alpha\,\alpha\,\alpha\,\bar{\alpha}\,\bar{\alpha}\,\alpha\,\bar{\alpha}\,\bar{\alpha}$ (Fig. 14(b)).

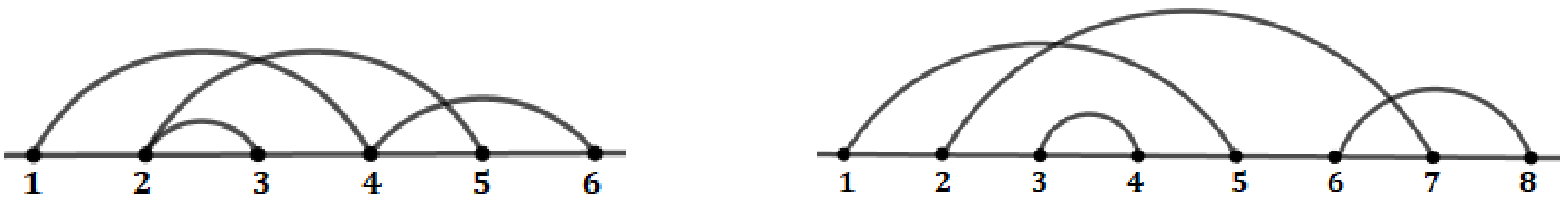

(a) (b)

Figure 14: $z = \alpha\, r\, \bar{\alpha}\, k\, \bar{\alpha}\, \bar{\alpha}$ $\quad a = \alpha\, \alpha\, \alpha\, \bar{\alpha}\, \bar{\alpha}\, \alpha\, \bar{\alpha}\, \bar{\alpha}$

**$k$-noncrossing**

We recall [6, p.24 ], that for any positive integer $k$, a $k$-tuple of distinct arcs $(i_1,j_1)$, $(i_2,j_2)$, ..., $(i_k,j_k)$ in a $b$-diagram (with $R=\bar{R}=K=\emptyset$ ) is a crossing $k$-tuple, iff $i_1<i_2< \ldots<i_k<j_1<j_2< \ldots<j_k$ .
A $b$-diagram is called $k$-noncrossing iff it contains non crossing $k$-tuple of arcs. For example, the $b$-diagram of Fig. 15 is a 4-noncrossing but not 3-noncrossing since it contains the crossing triple (3,6), (4,7) and (5,8).

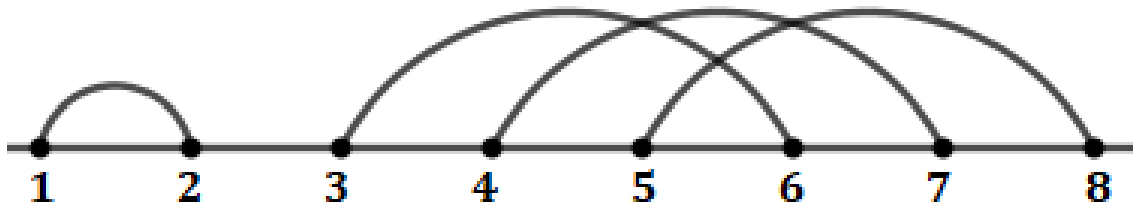

Figure 15: The 4-noncrossing b = 1 2 | 3 6 | 4 7 | 5 8

**Plato's transformations**

The transformations of Plato presented in his dialogue "Cratylus" 394A-B, i.e. addition, subtraction and transposition, which have successfully been applied from Alan Turing in his machine and from Noam Chomsky in Syntactic Structures, are also applied to $b$-diagrams cases.

## 4 Generation

A permutation $\sigma \in \Sigma_n$, is a generator of a $b$-diagram, when $b$ is obtained by removing a set of arcs $C \subset U_\sigma$ from the corresponding $\sigma$-diagram. For example, the permutation $\sigma = 1\ 3\ 2\ 7\ 8\ 4\ 5\ 6$ is a generator of the b-diagrams of Fig. 8(a) and Fig. 8(b).
It is obvious that the concatenation $b_1 b_2 \ldots b_m$ of the $m$ blocks of a $b$-diagram defines a permutation in $[n]$. The following proposition is valid:

<u>Proposition:</u> For every diagram $b = b_1 \mid b_2 \mid \ldots \mid b_m$ , the permutation $b_1 b_2 \ldots b_m$ in $[n]$, is one of its generators.

For example, $\sigma$=1 4 2 3 6 5 8 7 is the corresponding generator to the diagram $b$=1 4 | 2 | 3 6 | 5 8 | 7. Similarly, for the diagram $b$ = 2 1 3 4 | 5 8 7 6 , the permutation 2 1 3 4 5 8 7 6 gives the generator $\sigma$ = 1 3 4 5 6 7 6 2.

<u>Proposition:</u> For every diagram $b = b_1 \mid b_2 \mid \ldots \mid b_m$ , it is $|C| = n + m - \sum_{i=1}^{m} |b_i|$ .

<u>Proposition:</u> The number of generators of a $b$-diagram with $m$ blocks and $\ell$ isolated vertices equals to $2^{m-\ell}(m\text{-}1)!$

Proof: The number of generators coincides with the number of the possible concatenations which derive from the $m$ blocks (where the block $b_1$ is the first block) as well as of their reversal(1) ones, when they do not correspond to isolated vertices. From the first we get $(m\text{-}1)!$ generators and from the second $2^{m-\ell}$, therefore in total we get $2^{m-\ell}(m\text{-}1)!$ generators.
So, to the diagram $b$=2 3 1 4 | 5 8 7 6 where $m$=2 and $\ell$=0 correspond $2^2 \cdot (2\text{-}1)! = 4$ generators, i.e.: 1 3 2 5 8 7 6 4, 1 3 2 6 7 8 5 4 and their reverse ones 1 4 6 7 8 5 2 3 and 1 4 5 8 7 6 2 3 . Also, to the diagram $b$=1 4 | 2 | 3 6 | 5 8 | 7 where $m$=5 and $\ell$=2, correspond $2^3 \cdot 4! = 192$ generators.

---

*(1): The reversal of* $b_i = b_i^1 b_i^2 \ldots b_i^k$ *is* $\tilde{b}_i = b_i^k b_i^{k-1} \ldots b_i^1$.

By following the rationale of the previous propositions, we can find the set of the generators of a $b$-diagram. For the diagram $b$ =1 2 3│4 7 8│5 6 (Fig. 16) there exist $2^3 \cdot 2! = 16$ generators, half of which are the following permutations:

1 2 3 4 7 8 5 6
1 2 3 4 7 8 6 5
1 2 3 8 7 4 5 6
1 2 3 8 7 4 6 5
1 2 3 5 6 4 7 8
1 2 3 6 5 4 7 8
1 2 3 5 6 8 7 4
1 2 3 6 5 8 7 4

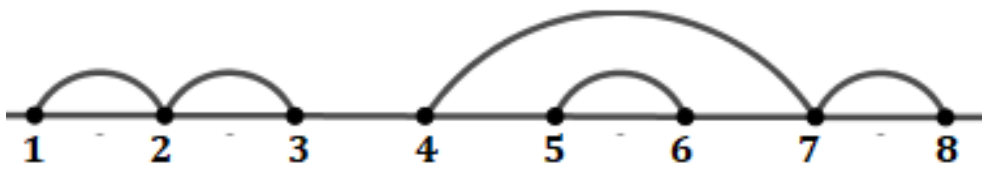

Figure 16: Diagram $b$= 1 2 3│4 7 8│5 6

and the other (8) half are their reverse ones.

Proposition: Two $b$-diagrams $b$ and $b'$ have common generators iff one of their respective sets $U_b$ and $U_{b'}$ is a subset of the other.

For example, for the diagrams of Fig. 17(a) $b$ =1 6│2 3│4 8 7│5 (which has $2^3 \cdot 3!=48$ generators) and of Fig. 17(b) $b$= 3 2 1 6│5 4 8 7 (which has $2^2 \cdot 1!=4$ generators) it is $U_b \subset U_{b'}$, it follows that they will have common generators, those of $b'$.

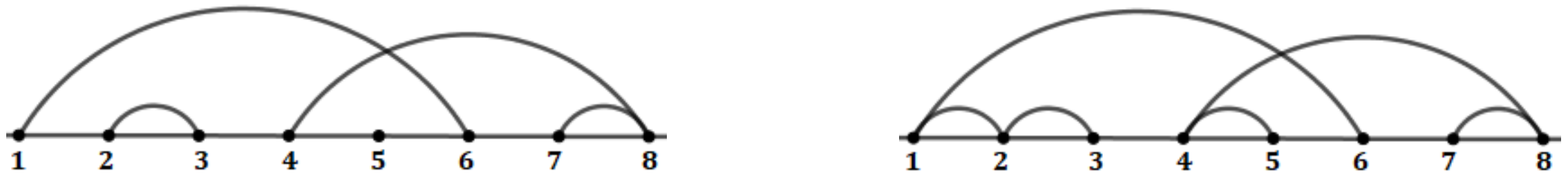

(a) (b)

Figure 17: $b$=1 6│2 3│4 8 7│5 $U_b$={16,23,48,5,78} $b$= 3 2 1 6│5 4 8 7 $U_{b'}$={12,16,23,45,48,78}

The identification of the generators of a $b$-diagram can also be done with the table of the ordered pairs of the respective set $U_b$ . So, to the diagram $b$= 1 2 3│4 7 8│5 6 of Fig. 18, corresponds the table:

| $i$ | 1 | 2 | 4 | 5 | 7 | • | • | • |
|---|---|---|---|---|---|---|---|---|
| $j$ | 2 | 3 | 7 | 6 | 8 | • | • | • |

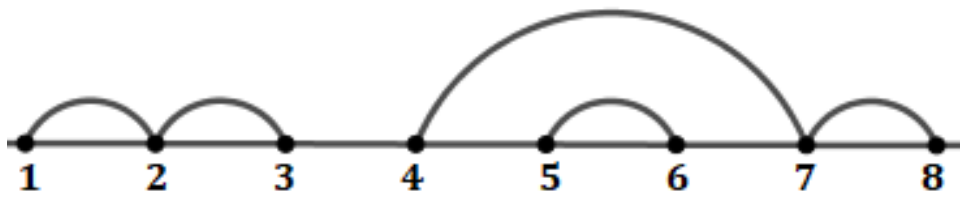

Figure 18: $b$= 1 2 3│4 7 8│5 6

It is obvious that each completion of the table's vacant cells defines a distinct generator of it. Given that 1 ϵ $R$ and 8 ϵ $\bar{R}$ , it is easy to derive that:

| $i$ | 1 | ② | 4 | 5 | ⑦ | 1 | • | • |
|---|---|---|---|---|---|---|---|---|
| $j$ | ② | 3 | ⑦ | 6 | 8 | • | • | 8 |

$R$= {1, . . . }, $\bar{R}$ = { . . . , 8} K = {2, 7}

So, in order to find the generators, it is each time necessary to complete only four cells with numbers from the set {3,4,5,6}. It is $i$ ϵ {{1,3,4}, {1,3,5}, {1,3,6}} because the number 3 already exists in the $j$ line, and respectively it is $j$ ϵ {{5,6,8}, {4,6,8}, {4,5,8}}. In the table below, the elements that define the completion of 8 generators (from the total $2^3 \cdot 2$=16), are denoted with an asterisk (*). The remaining 8 elements include cycles or previous results.

| $i$ | 1 | 3 | 4 | 1 | 3 | 5 | 1 | 3 | 6 |
|---|---|---|---|---|---|---|---|---|---|
| | 5 | 6 | 8 | 4 | 6 | 8* | 4 | 5 | 8 |
| | 5 | 8 | 6* | 4 | 8 | 6 | 4 | 8 | 5 |
| $j$ | 6 | 5 | 8 | 6 | 4 | 8* | 5 | 4 | 8 |
| | 6 | 8 | 5* | 6 | 8 | 4 | 5 | 8 | 4* |
| | 8 | 5 | 6* | 8 | 4 | 6 | 8 | 4 | 5 |
| | 8 | 6 | 5* | 8 | 6 | 4 | 8 | 5 | 4* |

Complement of a $b$-diagram in $[n]$ as to a generator $\sigma \in \Sigma_n$ of it, is called every $\bar{b}$-diagram in $[n]$ for which it is $U_{\bar{b}} = U_\sigma \setminus U_b$ . For example, for $b = 1\ 6\ 4 \mid 2\ 3\ 8 \mid 5\ 7$ (see Fig. 19(b)) and its generator $\sigma = 1\ 2\ 3\ 8\ 7\ 5\ 4\ 6$ (see Fig. 19(a)), the complement will be $\bar{b} = 1\ 2 \mid 3 \mid 4\ 5 \mid 6 \mid 78$ (see Fig. 19(c)).

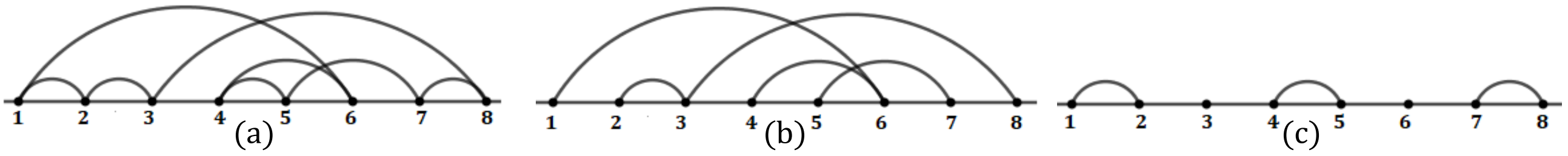

Fig. 19: generator $\sigma = 1\ 2\ 3\ 8\ 7\ 5\ 4\ 6$ $\quad$ $b$-diagram $b = 1\ 6\ 4 \mid 2\ 3\ 8 \mid 5\ 7$ $\quad$ complement $\bar{b} = 1\ 2 \mid 3 \mid 4\ 5 \mid 6 \mid 78$

$U_\sigma$={12,16,23,38,45,46,78} $\quad$ $U_b$={16,23,38,46,57} $\quad$ $U_{\bar{b}}$ = {12,3,45,6,78}

## Conclusions

The existence of a set of generators for every $b$-diagram enables the following:

– A new approach to the $b$-diagrams by using cyclic permutations and other tools of the Combinatorial science.
– The study of $b$-diagrams which have as generators the meandric colliers [4].
– The ability to identify relative $b$-diagrams and to classify them into equivalence classes.
– The study of mutation cases.
– New ways of presenting certain Biological methods.
– Alternative control methods of the credibility of biological results.

The Algebra of $b$-diagrams, combined with Plato's transformations, may provide viable solutions to certain biological problems, as to:

– The generation of new $b$-diagrams by the use of the transformations of addition and transposition.
– The classic problem of splicing which can be approached by adding arcs among different $b$-diagrams, thus creating connections between them.
– The removal of arcs from a $b$-diagram by using a vertical axis reminds of the biological scissors (crisptix).
– The kernel of a $b$-diagram is necessary for the study of the diagram's structure, since it defines its stability set as well as its absorbent set .
– The well-known acts of dilation and erosion in Topology [2], correspond fully to those of the determination of $\sigma$-diagrams and $b$-diagrams.
– The non-crossing problem acquires new dimensions with the $b$-diagrams.

**Acknowledgements** The author would like to thank Mr. Christos T. Kallinis for carefully reading the manuscript and providing many helpful comments and suggestions.